\documentclass[a4paper, 11pt,psamsfonts]{amsart}
\usepackage{amsmath}
\usepackage{amsthm}
\usepackage{amssymb}
\usepackage{amscd}
\usepackage{amsfonts}
\usepackage{amsbsy}
\usepackage{enumerate}
\usepackage{graphicx}
\usepackage{calc}
\usepackage{xcolor}
\usepackage{setspace}

\newtheorem {theorem} {Theorem}

\newtheorem {lemma}{Lemma}

\title[Local integrability of differential systems]
{A note on local integrability of differential systems\footnote{Published in J. Differential Equations}}

\author[X. Zhang]
{Xiang Zhang}

\address{School of Mathematical Sciences, MOE--LSC, Shanghai Jiao Tong University, Shanghai, 200240, P. R. China}
\email{xzhang@sjtu.edu.cn}

\subjclass[2010]{ 34A34; 34C20; 34C45; 37J30; 37K10.}

\keywords{Analytic differential systems; non--isolated singular point; local integrability; invariant manifold.}

\begin{document}

\begin{abstract}
For an $n$--dimensional local analytic differential system $\dot x=Ax+f(x)$ with $f(x)=O(|x|^2)$,
the Poincar\'e nonintegrability theorem states that if the eigenvalues of $A$ are not resonant, the system does not have an analytic or a formal first integral in a neighborhood of the origin. This result was extended in 2003 to the case when $A$ admits one zero eigenvalue and the other are non--resonant: for $n=2$ the system has an analytic first integral at the origin if and only if the origin is a non--isolated singular point; for $n>2$ the system has a formal first integral at the origin if and only if the origin is not an isolated singular point. However, the question of \emph{whether the system has an analytic first integral at the origin provided that the origin is not an isolated singular point} remains open.
\end{abstract}

\maketitle

\section{Introduction and statement of the main results}\label{s1}

For the local analytic differential system
\begin{equation}\label{e1}
\dot x=Ax+f(x), \qquad x\in \mathbb R^n,
\end{equation}
with $f(x)=O^*(\|x\|^2)\in C^\omega(\mathbb R^n,0)$, the study of the theory of local integrability or of the existence of first integrals at the origin can be traced back to Poincar\'e \cite{Po1891}.  Since then, the theory of local integrability has been greatly developed, see for example \cite{Bi1979, CYZ2008, CLZ2011, DRZ2016, DOR2016,Fu1996, LPW2012, RS2009,RXZ2014, SL2001, Sh2007,Zh2008, Zh2013}.
Hereafter, $O^*(\|x\|^2)$ denotes a function (or a vector--valued function) without constant and linear terms in its Taylor expansion, and $C^\omega(\mathbb R^n,0)$ denotes the set of analytic functions defined in a neighborhood of the origin.
Note that $x=0$ is a singular point of system \eqref{e1}, and that after an invertible linear change of coordinates we can always transform system \eqref{e1} to a system with its linear part matrix in Jordan normal form. So in what follows, we assume without loss of generality that $A$ is in Jordan normal form.

Let $\lambda=(\lambda_1,\ldots,\lambda_n)$ be the eigenvalues of the matrix $A$. Set
\[
\mathcal M_\lambda:=\left\{m\in\mathbb Z_+^n|\ \langle m,\,\lambda\rangle=0, \,\, |m|\ge 1\right\},
\]
where $\mathbb Z_+$ is the set of nonnegative integers, $\langle \cdot,\cdot\rangle$ denotes the inner product of two vectors in $\mathbb C^n$ and $|m|=m_1+\ldots+m_n$ for $m=(m_1,\ldots,m_n)$. If $\mathcal M_\lambda=\emptyset$, we say $\lambda$ is \textit{non--resonant}. If $\mathcal M_\lambda\ne\emptyset$,
each element of $\mathcal M_\lambda$ is called a \textit{resonant lattice}.

Poincar\'e \cite{Po1891} proved the next result.

\noindent\textbf{Theorem A.} {\it If system \eqref{e1} is analytic, and the eigenvalues $\lambda$ of $A$ is non--resonant, then the system has neither analytic nor formal first integrals.
}

For a proof of the Poincar\'e's result, see example \cite{Fu1996} or \cite{SL2001}.

Recall that a formal first integral is a formal series $H(x)$ which satisfies $\langle \nabla H(x),\, Ax+f(x)\rangle\equiv 0$ in $(\mathbb R^n,0)$, where $\nabla H$ is the gradient of $H$ and the partial derivative of $H$ is taken over all homogeneous terms in the sum of $H$, that is, if $H(x)=\sum\limits_{\ell=1}\limits^\infty H_\ell(x)$ with the $H_\ell$'s homogeneous polynomials of degree $\ell$ then $\nabla H=\sum\limits_{\ell=1}\limits^\infty \nabla H_\ell(x)$.

When the $n$--tuple of eigenvalues $\lambda$ are resonant, i.e. $\mathcal M_\lambda\ne\emptyset$, there are certain known results which provide necessary conditions ensuring the existence and number of functionally independent local analytic or formal first integrals of system \eqref{e1}. For more details, see \cite{CYZ2008,CLZ2011,DOR2016,LPW2012,RXZ2014}. About the equivalent characterization of analytic integrability via normal form, there are also some known results on the existence of analytic normalization of analytically integrable differential systems to their Poincar\'e--Dulac normal forms. For more details, see \cite{It1989, It2009, Zh2008, Zh2013, Zu2002, Zu2005}.
As our knowledge, there are very few general results providing necessary and sufficient conditions for the existence of analytic or formal first integrals defined in a neighborhood of the origin.

Li \textit{et al} \cite{LLZ2003} in 2003 studied the existence of local first integrals at the origin in case when one of the eigenvalues vanishes and the other are non--resonant, that is
\begin{equation}\label{e2}
\lambda_1=0 \ \mbox{ and } \  \sum\limits_{j=2}\limits^nm_j\lambda_j\ne 0\ \mbox{ for} \  m_j\in\mathbb Z_+  \mbox{ and} \ \sum\limits_{j=2}\limits^nm_j\ge 1.
\end{equation}
Their results can be stated as follows.

\noindent\textbf{Theorem B.} {\it Assume that the differential system \eqref{e1} is analytic and the conditions \eqref{e2} hold.
\begin{itemize}
\item[$(a)$] For $n>2$, system \eqref{e1} has a formal first integral in a neighborhood
of $x=0$ if and only if the singular point $x=0$ is not isolated. In particular,
if the singular point $x=0$ is isolated, system \eqref{e1} has no analytic first
integrals in a neighborhood of $x=0$.

\item[$(b)$] For $n=1,\,2$, system \eqref{e1} has an analytic first integral in a neighborhood of
$x=0$ if and only if the singular point $x=0$ is not isolated.
\end{itemize}
}

We note that statement $(b)$ of Li \textit{et al} \cite{LLZ2003} completely solve the problem on the existence of analytic first integral of a planar analytic differential system at its singular point with its linear part having eigenvalues satisfying \eqref{e2}. If system  \eqref{e1} is higher dimensional, the following problem remains open since 2003:

  {\it Whether is true that the analytic differential system \eqref{e1} under the condition \eqref{e2} has an analytic first integral in a neighborhood of the origin if and only if the singular point $x=0$ is not isolated?}

Here we give an answer to this problem. The main results are the following.

\begin{theorem}\label{t1}
For the local analytic differential system \eqref{e1}, assume that $\lambda_1=0$, and that $\lambda_2,\ldots,\lambda_n$ either all have positive real parts or all have negative real parts. Then system \eqref{e1} has an analytic first integral in $(\mathbb R^n,0)$ if and only if the singular point $x=0$ is not isolated.
\end{theorem}

Note that the conditions on the eigenvalues $\lambda$ in Theorem \ref{t1} is a special case of the conditions \eqref{e2}. Of course, the condition \eqref{e2} means that the eigenvalues $\lambda_2,\ldots$, $\lambda_n$ have no vanishing real parts. Statement $(b)$ of Theorem B is a special case of our Theorem \ref{t1}.

We remark that the method in \cite{LLZ2003} for proving the sufficiency of statement $(b)$ of Theorem B strongly depends on the dimension two, and it cannot be extended to higher dimensional systems. Here we provide a new approach using the ideas partly from \cite[Lemma 5.3]{Fe1979}.
The proof of the necessity of statement $(b)$ of Theorem B in \cite{LLZ2003} follows from that of statement $(a)$, which is much involved. Now we provide a new and simple proof using our newly developed theory of integrability given in \cite{DRZ2016}.

Next we present a $C^\infty$ version of Theorem \ref{t1}.

\begin{theorem}\label{t2}
For the $C^\infty$ differential system \eqref{e1}, assume that $\lambda_1=0$, and that $\lambda_2,\ldots,\lambda_n$ either all have positive real parts or all have negative real parts. Then system \eqref{e1} has a $C^\infty$ non--flat first integral in $(\mathbb R^n,0)$ if and only if there exists a curve passing through the origin which fulfils singular points of system \eqref{e1}.
\end{theorem}

Recall that a $C^\infty$ function is \textit{non--flat} if its Taylor expansion does not identically vanish. In $C^\infty$ category we cannot say, as stated in Theorem \ref{t1}, the singular point $x=0$ is not isolated, because a $C^\infty$ function of one variable can have zeros accumulating at one of its zeros, but it does not identically vanish. Whereas an analytic function of one variable either is identically equal to zero or has only isolated zeros.

Finally, we answer the open problem in the general condition \eqref{e2}, which shows that the sufficient part of the open problem is not correct in general.

\begin{theorem}\label{t2-1}
For the $C^r$ differential system \eqref{e1}, $r\in\{\omega,\infty\}$, suppose that the conditions \eqref{e2} hold.
\begin{itemize}
\item[$(a)$] If system \eqref{e1} has a $C^r$ first integral $($non--flat in the $C^\infty$ case$)$ in $(\mathbb R^n,0)$, then there exists a curve passing the origin, which fulfils singular points of system \eqref{e1}.

\item[$(b)$]  If system \eqref{e1} has a curve fulfilling singular points and containing the origin in its interior, the following statements hold.
\begin{itemize}
\item[$(b_1)$] System \eqref{e1} always has a formal first integral in $(\mathbb R^n,0)$.
\item[$(b_2)$] There exist analytic differential systems of the form \eqref{e1} which have no analytic first integrals in $(\mathbb R^n,0)$.
\end{itemize}
\end{itemize}
\end{theorem}

We remark that statement $(a)$ of Theorem \ref{t2-1} is in fact a different statement of Theorem B, here we add a little more information in it, because any analytic and $C^\infty$ function can be expanded as a Taylor series. Statement $(b_1)$ of Theorem \ref{t2-1} is the same as the sufficient part of Theorem B. Here we will present a different and simple proof, which can also be used to prove statement $(b_2)$. Statement $(b_2)$ of Theorem \ref{t2-1} provides a  negative answer to the sufficiency part of the open problem.

We do not know whether there exists a $C^\infty$ first integral if the conditions of statement $(b)$ hold.

This paper is organized as follows. In the next section we prove Theorem \ref{t1}. Sections \ref{s3} and \ref{s4} are the proofs of Theorems \ref{t2} and  \ref{t2-1}, respectively.

\section{Proof of Theorem \ref{t1}}\label{s2}

To prove this theorem we need the next result, due to Du {\it et al} \cite[Theorem 1 $(b)$]{DRZ2016}, which plays a central role in the proof of Theorem \ref{t1}.

\begin{theorem}\label{ta1}
Let $m\in\{1,\ldots,n-1\}$ be the maximum number of $\mathbb Q_+$--linearly independent elements of $\mathcal M_\lambda$. If the analytic differential system \eqref{e1} has $m$ functionally independent analytic or formal first integrals defined in a neighborhood of the origin, then it has $m$ functionally independent first integrals of the form
 \begin{equation}\label{e3}
    H_1(x)=x^{\alpha_1}+h_1(x),\,\ldots,\, H_{m}(x)=x^{\alpha_{m}}+h_{m}(x),
 \end{equation}
where $\alpha_1,\ldots,\alpha_{m}$ are $\mathbb Q_+$--linearly independent elements of $\mathcal M_\lambda$, and each $h_j(x)=O^*\left(|x|^{|\alpha_j|+1}\right)$, $j=1,\,\ldots,\, m$, consists of non--resonant monomials.
\end{theorem}

We now prove Theorem \ref{t1}. By the condition \eqref{e2} we can assume without loss of generality that system \eqref{e1} is of the form
\begin{equation}\label{e4}
\begin{array}{l}
\dot x_1=f_1(x),\\
\dot y=By+g(x),
\end{array}
\end{equation}
with $y=(x_2,\ldots,x_n)^T$ and $g(x)=(f_2(x),\ldots,f_n(x))^T$, where $T$ denotes the transpose of a vector or of a matrix.

\noindent{\it Necessity}. By assumption system \eqref{e1} has an analytic first integral in a neighborhood of the origin. The hypothesis on the eigenvalues $\lambda$ means that $\mathcal M_\lambda$ has a unique $\mathbb Q_+$--linearly independent element. So we are in the assumption of Theorem \ref{ta1}. It  follows that system \eqref{e1} has an analytic first integral of the form
\[
H(x)=x^\alpha+h(x),
\]
where $\alpha\in\mathcal M_\lambda$, and $h(x)$ is the higher order terms. Note that $(1,0,\ldots,0)$ is a basis of $\mathcal M_\lambda$. Take $\alpha=(1,0,\ldots,0)$, we obtain $H(x)=x_1+h(x)$.

Take the invertible change of coordinates $u=(u_1,\ldots,u_n)=\Phi(x)$ defined by
\[
u_1=H(x),\quad v=(u_2,\dots,u_n)^T:=(x_2,\ldots,x_n)^T,
\]
then system \eqref{e1} is transformed to
\begin{equation}\label{e5}
\begin{array}{l}
\dot u_1= \langle\nabla H,\, \dot x\rangle=\langle\nabla H,\, Ax+f(x)\rangle\equiv 0,\\
\, \dot v\,= \, Bv+g\circ \Phi^{-1}(u).
\end{array}
\end{equation}
Since $g=O^*(\|x\|^2)$ and $\Phi^{-1}(u)$ is near identity, it follows from the Implicit Function Theorem that system \eqref{e5} has singular points fulfilling the analytic curve defined by $Bv+g\circ \Phi^{-1}(u)=0$, which contains the origin as an interior point. By the invertible change of coordinates $u=\Phi(x)$, we get that system \eqref{e1} has an analytic curve passing through the origin, which fulfils singular points of system \eqref{e1}. This proves the necessity.

\noindent {\it Sufficiency}. By assumption it follows that the matrix $B$ in \eqref{e4} has no zero eigenvalues, so by the Implicit Functional Theorem it follows that the functional equation
\[
By+g(x)=0,
\]
with $x=(x_1,y)$, has an analytic solution, say $y=\varphi(x_1)$, defined in a neighborhood of $x_1=0$, which is tangent to the $x_1$--axis, where we have used the fact that $g(x)=O^*(\|x\|^2)$. Since by assumptions $f_1(x)$ is analytic and the singular points of system \eqref{e4} is not isolated, it forces that
\[
f_1(x_1,\varphi(x_1))\equiv 0.
\]
This means that under the sufficient assumption system \eqref{e1} has the analytic curve $y=\varphi(x)$ fulfilling singular points and including the origin in its interior.

Take the invertible and analytic change of coordinates
\begin{equation}\label{e6}
u_1=x_1,\quad v=(u_2,\ldots,u_n)^T:=y-\varphi(x_1)=(x_2,\ldots,x_n)^T-\varphi(x_1),
\end{equation}
and we denote by $u=(u_1,v)=G(x)$. Then system \eqref{e4} is transformed to the system
\begin{align}\label{e7}
\dot u_1=& \left\langle\int_0^1\partial yf_1(u_1,\alpha v+\varphi(u_1))d\alpha,\, v\right\rangle=R_1(u),\nonumber \\
\dot v=& Bv+\int_0^1\partial_yg(u_1,\alpha v+\varphi(u_1))d\alpha \, v\\
& \quad -\varphi'(u_1) \left\langle\int_0^1\partial yf_1(u_1,\alpha v+\varphi(u_1))d\alpha,\, v\right\rangle=R_2(u). \nonumber
\end{align}
Here we have used the next calculations that
\begin{align*}
f_1(x)=&f_1(u_1,v+\varphi(u_1))=\int_0^1\frac{d}{d\alpha}f_1(u_1,\alpha v+\varphi(u_1))d\alpha\\
 =& \left\langle\int_0^1\partial yf_1(u_1,\alpha v+\varphi(u_1))d\alpha,\, v\right\rangle,\\
\dot y-\varphi'(x_1)\dot x_1=&B(v+\varphi(u_1))+g(u_1,v+\varphi(u_1))-\varphi'(u_1)f_1(u_1,v+\varphi(u_1))\\
=&Bv+\int_0^1\frac{d}{d\alpha}\partial_yg(u_1,\alpha v+\varphi(u_1))d\alpha \,v\\
&\qquad\quad -\varphi'(u_1) \left\langle\int_0^1\partial yf_1(u_1,\alpha v+\varphi(u_1))d\alpha,\, v\right\rangle,
\end{align*}
where $\partial_yf_1$ and $\partial_yg$ denote the Jacobian matrices of $f_1(x,y)$ and $g(x,y)$ with respect to $y$, respectively. Moreover, we have
\begin{align*}
 &R_1(u)=O^*(\|u\|^2), \ \quad r_2(u):=R_2(u)-Bv=O^*(\|u\|^2), \\
 &R_1(u)|_{v=0}\equiv 0, \qquad \quad R_2(u)|_{v=0}\equiv 0.
\end{align*}

Obviously system \eqref{e7} is analytic, and it has the $u_1$--axis fulfilling singular points, which is called a  {\it singular line}. Since $A$ is real and $\lambda_2,\ldots,\lambda_n$ are non--resonant, it implies that the singular line $v=0$ is normally hyperbolic in a neighborhood of the origin.
Here \textit{normally hyperbolic} means that system \eqref{e7} has no eigenvalues with vanishing real parts in the directions normal to $v=0$.

For each singular point $(u_1,0)$ of system \eqref{e7} with $|u_1|$ suitably small, under the assumption of Theorem \ref{t1} the linearization of system \eqref{e7} at $(u_1,0)$ always has one zero eigenvalue and the other $n-1$ eigenvalues non--resonant. In fact, the $n-1$ eigenvalues either all
have positive real parts or all have negative real parts. Hence,
 by the Stable Manifold Theorem (see for example \cite[Theorem 3.2.1]{GH1983} and \cite[\S 4.1]{Ch2006}) it follows that system \eqref{e7} has an $(n-1)$--dimensional analytic stable (resp. unstable) invariant manifold at $(u_1,0)$ in case when all the eigenvalues of $B$ have negative (resp. positive) real parts, which is unique and tangent to the $(n-1)$--dimensional invariant linear space of the linearized system of system \eqref{e7} at $(u_1,0)$ associated to the eigenvalues with non--vanishing real parts. Moreover, it is clear that system \eqref{e7} has the unique and the same center manifold, i.e. $v=0$, at all singular points $(u_1,0)$. Furthermore, all points near the origin must belong to one of the $(n-1)$--dimensional invariant manifolds. This shows that the neighborhood of the origin is foliated by the family of $(n-1)$--dimensional analytic invariant manifolds, rooted at $(u_1,0)$ for $u_1\in (-\delta,\, \delta)$ for $\delta>0$.

In what follows, without loss of generality we consider only the case that all eigenvalues of $B$ have negative real parts. Denote by $u_1=\Phi_{u_1^0}(v)$ the $(n-1)$--dimensional invariant analytic manifold passing through $(u_1^0,0)$ with $u_1^0\in(-\sigma,\sigma)$. Take the coordinate change of variables
\begin{equation}\label{e8}
\xi_1=\Phi_{u_1}(v),\qquad \eta=(\xi_2,\ldots,\xi_n)^T=v.
\end{equation}
Note that $\bigcup\limits_{u_1\in(-\sigma,\sigma)}\Phi_{u_1}(v)$ is an $n$--dimensional subregion of $\mathbb C^n$ limited by the two analytic hypersurfaces
$\Phi_{-\sigma}(v)$ and $\Phi_{\sigma}(v)$, and it forms an $n$--dimensional analytic center--stable invariant manifold of system \eqref{e7} with the analytic center manifold $v=0$, i.e. the $u_1$--axis. This implies that $\Phi_{u_1}(v)$ is analytic not only in $v$ but also in $u_1$. Obviously, we get from the construction of $\Phi_{u_1}(v)$ that $\partial_{u_1}\Phi_{u_1}(v)$ is not zero for $u_1\in(-\sigma,\sigma)$ because two different hypersurfaces $\Phi_{u_1}(v)$'s do not intersect. These last proofs verify that the transformation \eqref{e8} is analytic and invertible in a neighborhood of the origin. We denote it by $\xi=W(u)$.

We claim that under the transformation \eqref{e8} system \eqref{e7} is changed to
\begin{equation}\label{e9}
\begin{array}{l}
\dot \xi_1= \, 0,\\
\dot \eta= \, B\eta+ q(\xi)\eta,
\end{array}
\end{equation}
where $\xi=(\xi_1,\eta)=(\xi_1,\xi_2,\ldots,\xi_n)$, and $q(\xi)$ is a matrix--valued analytic function of order $n-1$. Indeed,  the first equation in \eqref{e9} follows from the invariance of the hypersurface $\Phi_{u_1}(v)$ with each fixed $u_1$ under the flow of system \eqref{e7}, which can also be calculated as follows.
\begin{align*}
\dot \xi_1=&\partial_{u_1}\Phi_{u_1}(v)\dot u_1+\partial_v\Phi_{u_1}(v)\dot v\\
=&\partial_{u_1}\Phi_{u_1}(v)R_1(u)+\partial_v\Phi_{u_1}(v) R_2(u)\\
=&\langle \nabla \Phi_{u_1}(v), (R_1(u),R_2(u))\rangle=0,
\end{align*}
where $\nabla\Phi$ is the gradient of the function $\Phi$, and the last equality follows from the fact that the gradient $\nabla \Phi_{u_1}(v)$ is perpendicular to the tangent space of the invariant hypersurface at $u=(u_1,v)$ and the vector field $(R_1(u),R_2(u))$ belongs to the tangent space of the hypersurface at $u=(u_1,v)$. The second equation in \eqref{e9} follows easily from the expression of the second equation in \eqref{e7} and the fact $\eta=v$.

Clearly system \eqref{e9} has the first integral $V(\xi)=\xi_1$, it induces that system \eqref{e7} has the analytic first integral $V\circ W(u)$. Consequently system \eqref{e4} has the analytic first integral $V\circ W\circ G(x)$, so system \eqref{e1} has an analytic first integral. Here we have used the fact that if system $\dot x=\rho(x)$ has a smooth first integral $\omega(x)$, and it can be transformed to $\dot y=\sigma(y)$ via an invertible smooth transformation $x=\vartheta(y)$, then system $\dot y=\sigma(y)$ has the smooth first integral $\omega\circ \vartheta(y)$.

This completes the proof of the sufficient part and consequently completes the proof of Theorem 1.  \qed

\section{Proof of Theorem \ref{t2}}\label{s3}

Theorem \ref{t2} is very similar to Theorem \ref{t1}, but since the Taylor expansion of a $C^\infty$ function is not necessarily convergent, and even through it is convergent, its limit may not be equal to the given $C^\infty$ function, so we will use a different approach here than the analytic case.
For doing so, we need the next results.

The first one is on the normal forms of $C^\infty$ systems, see for instance \cite[Theorem 6.1]{Li2000} or \cite{IL1999}.

\begin{theorem}\label{t3}
Suppose that system \eqref{e1} is $C^\infty$. Let $\mu_1,\ldots,\mu_k,\mu_{k+1},\ldots,\mu_n$ be the eigenvalues of $A$, and satisfy $\mbox{\rm Re} \mu_j\ne 0$ for $j=1,\ldots, k$ and $\mbox{\rm Re} \mu_s=0$ for $s=k+1,\ldots,n$. If $\mu_1,\ldots,\mu_k$ are non--resonant, then for any $m\in\mathbb N$, system \eqref{e1} is $C^m$ equivalent to the system
\begin{equation}\label{e10}
\begin{array}{l}
\dot u=w(u),\\
\dot v=C(u) v+O^*(\|v\|^2),
\end{array}
\end{equation}
where $u\in\left(\mathbb R^{n-k},0\right)$, $v\in\left(\mathbb R^k,0\right)$, $w(0)=0$, $\partial_uw(0)$ has the eigenvalues $\mu_{k+1},\ldots,\mu_n$ and $C(0)$ has the eigenvalues $\mu_1,\ldots,\mu_k$.
\end{theorem}

Recall by definition that \textit{two systems of the form \eqref{e1} are $C^m$ equivalent} if there exists a near identity $C^m$ transformation which sends one system to another. A \textit{near identity transformation} is the one of the form $x=(u,v)+O^*(\|(u,v)\|^2)$.

The second one is on the spectrum of a linear differential operator, see for example \cite[Lemma 1.1]{Bi1979} or \cite[Lemma 4.5]{Li2000}.

\begin{lemma}\label{l1}
Let $\mathcal H^r_n(\mathbb C)$ be the linear space formed by homogeneous polynomials of degree $r$ in $n$ variables with coefficients in $\mathbb C$. For two $n$th order matrices $M_1$ and $M_2$, we define a linear operator on $\mathcal H^r_n(\mathbb C)$ by
\[
\mathcal L(h)(x) =\left\langle \nabla h(x),M_1x\right\rangle-M_2h(x),\qquad h\in \mathcal H^r_n(\mathbb C).
\]
Then the spectrum of $\mathcal L$ on $\mathcal H^r_n(\mathbb C)$ is
\[
\sigma(\mathcal L):=\{\langle k,\mu\rangle-\nu_j|\ k\in \mathbb Z_+^n,\ |k|=r, \, j\in\{1,\ldots,n\}\},
\]
where $\mu=(\mu_1,\dots,\mu_n)$ and $\nu=(\nu_1,\ldots,\nu_n)$ are respectively the $n$--tuples of eigenvalues of $M_1$ and $M_2$.
\end{lemma}

We now prove Theorem \ref{t2}. As in the proof of Theorem \ref{t1} we assume without loss of generality that system \eqref{e1} has the form \eqref{e4}.

\noindent{\it Necessity}. By assumption let $H(x)$ be a $C^\infty$ non--flat first integral of system \eqref{e4}. We choose $m\in\mathbb N$ sufficiently large such that the Taylor expansion of $H(x)$ is of the form $H(x)=H_\ell(x)+O^*(\|x\|^{\ell+1})$, $\ell < m$, with $H_\ell(x)$ a homogenous polynomoal of degree $\ell$, which does not identically vanish.

Under the assumption of Theorem \ref{t2} we have only one zero eigenvalue, and the other eigenvalues are non--resonant. Hence, by Theorem \ref{t3} there exists a near identity $C^m$ transformation, say $x=\Phi(u)$, which sends system \eqref{e4} to a system of the form \eqref{e10}, i.e.
\begin{equation}\label{e11}
\dot u_1=w(u_1), \qquad
\dot v=C(u_1)v+O^*(\|v\|^2),
\end{equation}
where $u=(u_1,v)$, $v=(u_2,\ldots,u_n)$, $w(u_1)=O^*(u_1^2)$ and $C(0)$ has the eigenvalues $\lambda_2,\ldots,\lambda_n$. Then system \eqref{e11} has the $C^m$ first integral $\widetilde H(u):=H\circ\Phi(u)=H_\ell(u)+O^*(\|u\|^{\ell+1})$, where we have used the fact that $\Phi(u)$ is near identity.

By the properties of first integrals we get
\begin{equation}\label{e12}
\left\langle \partial_u\widetilde H(u),\, \left(w(u_1), C(u_1) v+O^*(\|v\|^2) \right)\right\rangle\equiv 0\quad \mbox{\rm in } \left(\mathbb R^n,0\right).
\end{equation}
Comparing the terms of degree $\ell$ in \eqref{e12} gives
\begin{equation}\label{e13}
\left\langle \partial_u H_\ell(u),\, \left(0,C(0) v\right)\right\rangle\equiv 0.
\end{equation}
Since the eigenvalues $\lambda_2,\ldots,\lambda_n$ of $C(0)$ are non--resonant, it follows from Lemma \ref{l1} that equation \eqref{e13} has only the solution of the form $H_\ell(u)=H_\ell(u_1)$. Moreover, the homogeneity of $H_\ell(u)$ forces that $H_{\ell}(u_1)=a_\ell u_1^\ell$ with $a_\ell \ne 0$ a constant.

Set $w(u_1)=b u_1^\sigma+O^*(u_1^{\sigma+1})$ with $\sigma\in\mathbb N\setminus\{1\}$ and $b$ a constant. Balancing the coefficients of the terms of the lowest degree in $u_1$ in $O(v^0)$ in \eqref{e12}, one gets that $a_\ell b=0$. This induces that $b=0$, and consequently $w(u_1)\equiv 0$.

The above proof shows that the line $v=0$ fulfils singular points of system \eqref{e11}. Hence, system \eqref{e4} has singular points fulfilling the curve $v\circ \Phi^{-1}(x)=0$, where $\Phi^{-1}(x)$ is the inverse of the transformation sending system \eqref{e4} to system \eqref{e11}. This proves the necessity.

\noindent{\it Sufficiency}.  The Implicit Function Theorem verifies that $By+g(x)=0$ has a unique solution, say $y=\varphi(x_1)$, which is $C^\infty$. Moreover, we have from the assumption that $f_1(x_1,\varphi(x_1))\equiv 0$.

Applying the Stable Manifold Theorem to system \eqref{e4} at each singular point $(x_1,\varphi(x_1))$ with $x_1\in(-\delta,\delta)$ for some small positive $\delta$, one obtains a unique and $(n-1)$--dimensional $C^\infty$ stable (resp. unstable) manifold depending on all $\mbox{\rm Re}\lambda_j>0$ (resp. all  $\mbox{\rm Re}\lambda_j<0$) for $j=2,\ldots,n$. We denote this invariant manifold by $\mathcal M_{x_1}(y)$. Then working in a similar way as in the proof of Theorem \ref{t1} and using the fact that the center manifold, i.e. $y=\varphi(x)$, is one dimensional and $C^\infty$,  we obtain that system \eqref{e1} has a $C^\infty$ first integral in a neighborhood of the origin. Furthermore, the $C^\infty$ first integral is not flat because the transformation, which is given in the proof of Theorem \ref{t1}, is not flat.

This completes the proof of Theorem \ref{t2}. \qed

\section{Proof of Theorem \ref{t2-1}}\label{s4}

\noindent{\it Proof of statement $(a)$}. This statement can be proved using the same arguments as those in the proof of Theorem \ref{t1} when the system is analytic, and of Theorem \ref{t2} when the system is $C^\infty$. It can also be obtained as a consequence of the necessary part of statement $(a)$ of Theorem B. The details are omitted.

\noindent{\it Proof of statement $(b)$}.  Statement $(b_1)$ is just the sufficient part of statement $(a)$ of Theorem B, which is well known. Here we provide a different and simple proof to it, which will be used in the proof of statement $(b_2)$.

Using the same arguments as in the proof of Theorem \ref{t1}, we can assume without loss of generality that system \eqref{e1} has the form \eqref{e7} with $B$ in Jordan normal form, and
\begin{equation}\label{e3-0}
f_1(x)|_{y=0}=0\  \mbox{ and } \ g(x)|_{y=0}=0.
\end{equation}
Here we have used the fact that systems \eqref{e1} and \eqref{e7} are $C^r$ equivalent for $r\in\{\infty,\omega\}$, and so they either both have a formal first integral or both have no a formal first integral. Hence we can write system \eqref{e1} in the form \eqref{e7} replacing $u$ and $v$ by $x$ and $y$, respectively.

Let $H(x)$ be an analytic function or a formal series with the expression
\begin{equation}\label{e3-1}
H(x)=\sum\limits_{j=m}\limits^\infty H_j(x),
\end{equation}
where $m\in\mathbb N$ and the $H_j(x)$ is a homogeneous polynomial of degree $j$. In order for $H(x)$ to be a first integral, we must have
\begin{equation}\label{e3-2}
f_1(x)\frac{\partial H}{\partial x_1} +\left\langle By+g(x),\, \frac{\partial H}{\partial y}\right\rangle\equiv 0.
\end{equation}
Let the Taylor expansions of $f_1(x)$ and $g(x)$ be
\begin{equation}\label{e3-3}
f_1(x)=\sum\limits_{s=2}\limits^\infty f_{1s}(x), \qquad
g(x)=\sum\limits_{s=2}\limits^\infty g_{s}(x),
\end{equation}
with the $f_{1s}$ homogeneous polynomial of degree $s$ and the $g_s(x)$ $(n-1)$--dimensional vector--valued homogeneous polynomial of degree $s$.

Substituting \eqref{e3-1} and \eqref{e3-3} into \eqref{e3-2}, and comparing the terms which have the same degree, one gets that
\begin{align}\label{e3-4}
\mathcal L^*(H_m)=&0,\\
\mathcal L^*(H_{m+\ell})=&
-\sum\limits_{j=2}\limits^{\ell+1}\left\langle g_j(x),\frac{\partial H_{m+\ell+1-j}}{\partial y}\right\rangle\nonumber\\
&\qquad\  -\sum\limits_{j=2}\limits^{\ell+1}  f_{1j}(x) \frac{\partial H_{m+\ell+1-j}}{\partial x_1},\quad \ell=1,2,\ldots\label{e3-5}
\end{align}
where $\mathcal L^*$ is the linear operator defined by
\[
\mathcal L^*=\left\langle By,\, \frac{\partial }{\partial y}\right\rangle.
\]
According to Lemma \ref{l1} the spectrum of $\mathcal L^*$ on $\mathcal H_{n-1}(y)$, the set of homogeneous polynomials in the $n-1$ variables $y$, is
\begin{equation}\label{e3-6}
\{\langle\lambda^*,\, m^*\rangle|\ m^*\in\mathbb Z_+^{n-1}, \, |m^*|\ge 1\},
\end{equation}
where $\lambda^*=(\lambda_2,\ldots,\lambda_n)$ and $m^*=(m_2,\ldots,m_n)$.

Since the eigenvalues $\lambda_2,\ldots,\lambda_n$ are non--resonant, it follows from the spectrum \eqref{e3-6} of $\mathcal L^*$ that equation \eqref{e3-4} has only the solution of the form $H_m(x)=H_m(x_1)$. Since $H_m(x)$ is homogeneous, it must be of the form $H_m(x)=a_mx_1^m$ with $a_m$ a non--zero constant.
For $\ell=1$, equation \eqref{e3-5} is reduced to
\[
\mathcal L^*(H_{m+1})=
  -   f_{12}(x) \frac{\partial H_{m}(x_1)}{\partial x_1}.
\]
Since $f_{12}(x)|_{y=0}\equiv 0$, by Lemma \ref{l1} this last equation has a unique solution $H_{m+1}(x)$ modulo a monomial $a_{m+1}x^{m+1}$ with $a_{m+1}$ an arbitrary constant.

To apply the induction, we assume that for $\ell =1,\ldots, k-1$, equation \eqref{e3-5} has a homogeneous polynomial solution $H_{m+\ell}(x)$, which is successively uniquely determined modulo a monomial $a_{m+\ell}x^{m+\ell}$ with $a_{m+\ell}$ a constant. For $\ell =k$,
we get from \eqref{e3-0} that
\[
g_j(x)|_{y=0}\equiv 0,\quad f_{1j}(x)|_{y=0}\equiv 0, \quad j=2,3,\ldots,k+1.
\]
This means that the right--hand side of equation \eqref{e3-5} identically vanishes when $y=0$, that is, each monomial in the right--hand side of equation \eqref{e3-5} is of the form $c_{p,q}x_1^py^q$ with $p\in\mathbb Z_+$, $q=(q_2,\ldots,q_n)\in\mathbb Z_+^{n-1}$ and $|q|\ge 1$. So, we obtain from the spectrum \eqref{e3-6} of $\mathcal L^*$ via Lemma \ref{l1} that equation \eqref{e3-5} with $\ell=k$ has a unique solution $H_{m+k}(x)$ modulo a monomial $a_{m+k}x^{m+k}$ with $a_{m+k}$ a constant.

By induction, for all $\ell\in\mathbb N$ equation \eqref{e3-5} has a unique solution $H_{m+\ell}(x)$ modulo a monomial $a_{m+\ell}x^{m+\ell}$ with $a_{m+\ell}$ a constant. This proves that system \eqref{e7} has a formal first integral in a neighborhood of the origin, and consequently statement $(b_1)$ follows.

To prove statement $(b_2)$, we choose the non--resonant eigenvalues $\lambda^*:=(\lambda_2,\ldots,\lambda_n)$ of $A$ such that the matrix $B$ in \eqref{e4} is diagonal and the set
\[
\mathcal D:=\{\langle \lambda^*,m^*\rangle|\ m^*\in\mathbb Z_+^{n-1},\, |m^*|\ge 1\}
\]
has the accumulation point $0$, and the points in a subset of $\mathcal D$ accumulates $0$ in extremely fast speed.

Note from the proof of $(b_1)$ that the solution $H_{m+\ell}(x)$ of the equation \eqref{e3-5} consists of $a_{m+\ell}x_1^{m+\ell}$ plus the sum of the monomials in the right--hand side of \eqref{e3-5} multiplied by a factor of the form $\langle \lambda^*,m^*\rangle^{-1}$. Besides, the $f(x)=(f_1(x),g(x))=O^*(\|x\|^2)$ in \eqref{e4} can be any $n$--dimensional vector--valued analytic function which satisfies $f(x)|_{y=0}\equiv 0$. This implies that we can choose analytic functions $f(x)$ such that the monomials $x_1^py^{m^*}$ appearing in the right--hand side of equation \eqref{e3-6} have their exponents $m^*$ satisfying that the set $\{\langle \lambda^*,m^*\rangle\}$ accumulates $0$ too rapidly, which leads to the series $\sum\limits_{s=m}\limits^\infty H_s(x)=\sum\limits_{s=m}\limits^\infty H_s(x_1,y)$ does not convergent in any small neighborhood of the origin. Consequently the statement follows.

For precise, we present a concrete example. Consider the case $n=3$ with $x=(x_1,y)\in\mathbb R^3$ and system \eqref{e4} satisfying that
\begin{itemize}
\item $f_1(x)=f_1(y)$ is an analytic function whose Taylor expansion contains all possible monomials of degree greater than $1$;
\item $B=\mbox{\rm diag}(1,-\zeta)$ with $\zeta$ a Liouville number;
\item $g(x)\equiv 0$.
\end{itemize}
Recall that a \textit{Liouville number} is an irrational number $c$ with the property that, for every $m\in\mathbb N$, there exist positive integers $p$ and $q$ with $q > 1$ such that
\[
\displaystyle 0<\left|c-{\frac {p}{q}}\right|<{\frac {1}{q^{m}}}.
\]
Since $g(x)\equiv 0$, equation \eqref{e3-5} is reduced to
\[
\mathcal L^*(H_{m+\ell})=
 -\sum\limits_{j=2}\limits^{\ell+1}  f_{1j}(y) \frac{\partial H_{m+\ell+1-j}(x)}{\partial x_1},\quad \ell=1,2,\ldots
\]
Since $\mathcal L^*$ is invertible on $\mathcal H_2(y)$, it follows that
\[
H_{m+\ell}(x)=H^{(1)}_{m+\ell}(x)+H^{(2)}_{m+\ell}(x),
\]
with
\begin{align*}
H^{(1)}_{m+\ell}(x)&=-\sum\limits_{j=2}\limits^{\ell} \left(\mathcal L^*\right)^{-1}\left( f_{1j}(y) \frac{\partial H_{m+\ell+1-j}(x)}{\partial x_1}\right)+a_{m+\ell} x_1^{m+\ell},\\
H^{(2)}_{m+\ell}(x)&=-m a_m x_1^{m-1}\sum\limits_{|m^*|=\ell+1} \frac{a_{m^*}^{(\ell+1)}}{m_2-m_3\zeta} y^{m^*},
\end{align*}
where $m^*=(m_2,m_3)$, and the $a_{m^*}^{(\ell+1)}$'s are the coefficients of the monomials $y^{m^*}$ in the homogeneous polynomial $f_{1,\ell+1}(y)$ and satisfy
\[
\left|a_{m^*}^{(\ell+1)}\right|=\dfrac{1}{\ell+1}\left(\dfrac{3}{2}\right)^{-(\ell+1)}.
\]

We claim that the formal first integral $H(x)$ does not converge in any neighborhood of the origin. On the contrary, we assume that $H(x)$ is convergent in $\|x\|<\rho$ with $\rho>0$. Then for any $\sigma\in(0,\rho)$ the formal first integral $H(x)$ is absolutely and uniformly convergent in $\|x\|\le\sigma$. This forces that the series
\[
H^{(2)}(x):=\sum\limits_{s=m+1}\limits^\infty H^{(2)}_s(x)=-m a_m x_1^{m-1}\sum\limits_{\ell=1}\limits^\infty \left(\sum\limits_{|m^*|=\ell+1} \frac{a_{m^*}^{(\ell+1)}}{m_2-m_3\zeta} y^{m^*}  \right)
\]
is convergent in  $\|x\|\le\sigma$. But it is impossible because $\zeta$ is a Liouville number, and consequently the set $\{m_2-m_3\zeta|\ m^*=(m_2,m_3)\in\mathbb Z_+^2, \, |m^*|\ge 2\}$ has a subset which accumulates $0$ in at least the speed $2^{-k}$ for $k\in\mathbb N$ and $k\rightarrow \infty$. This contradiction implies that the claim holds.

This last claim verifies that system \eqref{e4} under the given conditions has no analytic first integrals in a neighborhood of the origin.

This completes the proof of Theorem \ref{t2-1}. \qed

\section*{Acknowledgements}

We thank the referee for his/her nice comments which greatly improve the presentation of our paper.

The author is partially supported by NNSF of China grant
numbers 11271252 and 11671254, and by Innovation Program of Shanghai Municipal Education Commission
grant 15ZZ012.

\end{document}